\documentclass[12pt]{article}
\usepackage[dvips]{graphicx}
\DeclareGraphicsExtensions{ps,eps}
\makeatletter
\def\section{\@startsection{section}{1}{\z@}{-1.5ex plus -.5ex
minus -2.ex}{1ex plus .2ex}{\large\bf}}
{

\input{amssym.def}
\input{amssym}

\def\rit{{\Bbb R}}

\def\nit{{\Bbb N}}
\def\zit{{\Bbb Z}}
\def\tit{{\Bbb T}}

\def\@thmcounterstep{}
\long\def\@makecaption#1#2{\vskip 10pt \setbox\@tempboxa\hbox{#1.#2}
\ifdim \wd\@tempboxa >\hsize
#1.#2\par 
\else99
\hbox to\hsize{\hfil\box\@tempboxa\hfil}
\fi}
\def\ps@headings{
\def\@oddhead{\footnotesize\rm\hfill\runninghead\hfill}
\def\@evenhead{\@oddhead}
\def\@oddfoot{\rm\hfill\thepage\hfill}\def\@evenfoot{\@oddfoot}}
\newtheorem{Theorem}{Theorem}[section]
\newtheorem{Lemma}[Theorem]{Lemma}
\newtheorem{Definition}[Theorem]{Definition}
\newtheorem{Proposition}[Theorem]{Proposition}

\makeatother
\title
{
$L^2$ formulation of multidimensional scalar conservation laws
}
\def\runninghead{\quad 
$L^2$  formulation of scalar conservation laws
}
\author{
{\em Yann Brenier}\thanks{CNRS, Universit\'e de Nice, FR 2800,
brenier@math.unice.fr}}
\date{} 
\begin{document}
\pagestyle{headings}
\flushbottom
\maketitle
\vspace{-10pt}

\subsection*{}
Keywords : scalar conservation laws, level sets,
kinetic approximation,  maximal monotone operator
\\
\\
AMS classification~: 35L65, 47H05

\section*{Abstract}

We show that Kruzhkov's theory of entropy solutions to multidimensional
scalar conservation laws \cite{Kr} can be entirely recast in $L^2$
and fits into the general theory of maximal 
monotone operators in Hilbert spaces. 
Our approach is based on a
combination of level-set, kinetic and transport-collapse
approximations, in the spirit of previous works by
Giga, Miyakawa, Osher, Tsai and the author \cite{Br1,Br2,Br3,Br4,GM,TGO}.

\section{A short review of Kruzhkov's theory}

First order systems of conservation laws read:
$$
\partial_t u+\sum_{i=1}^d \partial_{x_i}(Q_i(u))=0,
$$
or, in short, using the nabla notation,
\begin{equation}
\label{scalar}  
\partial_t u+\nabla_x\cdot (Q(u))=0,
\end{equation}
where $u=u(t,x)\in \rit^m$ depends on $t\ge 0$, $x\in \rit^d$,
and $\cdot$ denotes the inner product in $\rit^d$.
The $Q_i$ (for $i=1,\cdot\cdot\cdot,d$) are given smooth functions
from $\rit^m$ into itself. The system is called hyperbolic when,
for each $\tau\in\rit^d$ and each $U\in\rit^m$, the $m\times m$ matrix
$\sum_{i=1,d}\tau_i Q'_i(U)$ can be put in diagonal
form with real eigenvalues.
There is no general theory to solve globally in time the initial value problem
for such systems of PDEs.
(See \cite{BDLL,Da,Ma,Se} for a general introduction to the field.)
In general, smooth solutions are known to exist
for short times but are expected to blow up in finite
time. Therefore, it is usual to consider discontinuous weak solutions,
satisfying additional 'entropy' conditions,
to adress the initial value problem, but nothing is
known, in general, about their existence . Some special situations
are far better understood. First, for some special systems 
(enjoying 'linear degeneracy'
or 'null conditions'), smooth solutions may be global (shock free), 
at least for 'small' initial data (see \cite{Kl}, for instance).
Next, in one space dimension $d=1$, for a large class of systems, 
existence and uniqueness of global weak entropy solutions have been 
(recently) proven for initial data
of sufficiently small total variation \cite{BB}.
Still, in one space dimension, for a limited class of systems (typically
for $m=2$), existence of global weak entropy solutions have been obtained
for large initial data by 'compensated compactness' arguments
\cite{Ta,Di,LPS}.
Finally, there is a very comprehensive theory in the much simpler
case of a $single$ conservation laws, i.e. when $m=1$.
Then, equation (\ref{scalar}) is called a 'scalar
conservation law'. Kruzhkov \cite{Kr} showed that 
such a scalar conservation law has a unique 'entropy
solution' $u\in L^\infty$ for each given initial condition $u_0\in L^\infty$.
(If the derivative $Q'$ is further
assumed to be bounded, then we can substitute $L^1_{loc}$ for $L^\infty$ in
this statement.)
An entropy (or Kruzhkov) solution is an $L^\infty$ function that satisfies the
following distributional inequality
\begin{equation}
\label{entropy}  
\partial_t C(u)+\nabla_x\cdot (Q^C(u))\le 0,
\end{equation}
for all Lipschitz convex function $C:\rit\rightarrow \rit$, 
where the derivative of $Q^C$ is defined by $(Q^C)'=C'Q'$.
In addition, the initial condition $u_0$ is prescribed
in $L^1_{loc}$, namely:
\begin{equation}
\label{continuity}  
\lim _{t\rightarrow 0} \int_{B} |u(t,x)-u_0(x)|dx=0,
\end{equation}
for all compact subset $B$ of $\rit^d$.
Beyond their existence and uniqueness,
the Kruzhkov solutions enjoy many interesting properties. Each entropy
solution $u(t,\cdot)$, with initial condition $u_0$, 
continuously depends on $t\ge 0$ in $L^1_{loc}$ and
can be written $T(t)u_0$,
where $(T(t),t\ge 0)$ is a family of order preserving operators:
\begin{equation}
\label{order T}
T(t)u_0\;\ge\; T(t)\tilde u_0\;,\;\;\forall t\ge 0,
\end{equation}
whenever $u_0\ge \tilde u_0$. Since constants are trivial
entropy solutions to (\ref{scalar}), it follows that 
if $u_0$ takes its values in
some fixed compact interval, so does $u(t,\cdot)$ for all $t\ge 0$.
Next, two solutions $u$ and $\tilde u$, with
$u_0-\tilde u_0\in L^1$, are $L^1$ stable with
respect to their initial conditions:
\begin{equation}
\label{L1 stability}
\int|u(t,x)-\tilde u(t,x)|dx
\le
\int|u_0(x)-\tilde u_0(x)|dx,
\end{equation}
for all $t\ge 0$.
As a consequence, the total variation $TV(u(t,\cdot))$
of a Kruzhkov solution $u$ at time $t\ge 0$ cannot be larger than the total
variation of its initial condition $u_0$.
This easily comes from the translation invariance of
(\ref{scalar}) and from the following definition of the total variation 
of a function $v$:
\begin{equation}
\label{TV}
TV(v)=\sup_{\eta\in\rit^d,\;\;\eta\ne 0}
\int\frac{|v(x+\eta)-v(x)|}{||\eta||}dx,
\end{equation}
where $||\cdot||$ denotes the Euclidean norm on $\rit^d$.
The space $L^1$ plays a key role in Kruzhkov's theory. 
There is no $L^p$ stability with respect to initial conditions in any $p>1$.
Typically, for $p>1$,
the Sobolev norm $||u(t,\cdot)||_{W^{1,p}}$ of a Kruzhkov
solution blows up in finite time.
This fact has induced a great amount of pessimism about the possibility
of a unified theory of global solutions for general
multidimensional systems of hyperbolic conservation laws. 
Indeed, simple linear
systems, such as the wave equation (written as a first order system)
or the Maxwell equations, are not well posed in any $L^p$ but for $p=2$
\cite{Brn}. However, as shown in the present work, $L^2$ is a perfectly
suitable space for entropy solutions to multidimensional scalar conservation
laws, provided a different formulation
is used, based on a
combination of level-set, kinetic and transport-collapse
approximations, in the spirit of previous works by
Giga, Miyakawa, Osher, Tsai and the author \cite{Br1,Br2,Br3,Br4,GM,TGO}.

\section{Kruzhkov solutions revisited}
\subsection{A maximal monotone operator in $L^2$}

Subsequently, we restrict ourself, for simplicity, to initial
conditions $u_0(x)$ valued in $[0,1]$ and spatially periodic of
period 1 in each direction. In other words, the variable $x$ will
be valued in the flat torus $\tit^d=\rit^d/\zit^d$.
\\
\\
Let us now introduce:
\\
1) the space $L^2([0,1]\times \tit^d)$ of all square
integrable functions
$$(a,x)\in [0,1]\times \tit^d\rightarrow Y(a,x)\in\rit\;,$$
2) the closed convex cone $K$ of all $Y\in L^2$ such that
$\partial_a Y\ge 0$ (in the sense of distributions),
\\
3) the subdifferential of $K$ defined at each point
$Y\in K$ by:
\begin{equation}
\label{subdifferential}
\partial K(Y)=\{Z\in L^2,\;\;\;
\int (\tilde Y-Y)Z\;dadx\le 0\;,
\;\;\;\forall \tilde Y\in K\}\;,
\end{equation}
4) the maximal monotone operator (MMO) (see \cite{Brz}):
\begin{equation}
\label{operator}
Y\rightarrow - q(a)\cdot\nabla_x Y+\partial K(Y),
\end{equation}
where $q(a)=Q'(a)$, and the corresponding
subdifferential equation \cite{Brz}:
\begin{equation}
\label{inclusion}
0\;\in\;\;
\partial_t Y+ q(a)\cdot\nabla_x Y+\partial K(Y).
\end{equation}

From maximal monotone operator theory \cite{Brz}, we know that,
for each initial condition $Y_0\in K$, there is a unique solution 
$Y(t,\cdot)\in K$ to (\ref{inclusion}), for all $t\ge 0$.
More precisely, we will use the following definition (which includes the possibility
of a left-hand side $q_0\in L^2([0,1])$):

\begin{Definition}
\label{def}
$Y$ is a solution to 
\begin{equation}
\label{inclusion bis}
q_0(a)\;\in\;\;
\partial_t Y+ q(a)\cdot\nabla_x Y+\partial K(Y),
\end{equation}
with initial value $Y_0\in K$ and left-hand side $q_0\in L^2([0,1])$,
if:
\\
1) $t\rightarrow Y(t,\cdot)\in L^2$ is continuous 
and valued in $K$,
with $Y(0,\cdot)=Y_0,$
\\
2) $Y$ satisfies, in the sense of distribution,
\begin{equation}
\label{semi-integral}
\frac{d}{dt} \int |Y-Z|^2 dadx
\le
2\;\int (Y-Z)(q_0(a)-\partial_t Z- q(a)\cdot\nabla_x Z)dadx,
\end{equation}
for each smooth function $Z(t,a,x)$
such that $\partial_a Z\ge 0$.

\end{Definition}

\begin{Proposition}
\label{Proposition}
For each $Y_0\in K$, and $q_0\in L^2([0,1])$, there is a unique solution $Y$ 
to (\ref{inclusion bis}) in the sense of Definition \ref{def}.
If both $Y_0$ and $q_0$ belong to $L^\infty$, then we have for all $t\ge0$:
\begin{equation}
\label{maxi}
-t\sup (-q_0)_++\inf Y_0\le Y(t,\cdot)\le \sup Y_0+t\sup (q_0)_+.
\end{equation}
If $\nabla_x Y_0$ belongs to $L^2$, 
then so do $\partial_t Y(t,\cdot)$ and $\nabla_x Y(t,\cdot)$ 
for all $t\ge 0$.
Two solutions $Y$ and $\tilde Y$ to (\ref{inclusion bis}) (with different
left-hand side $q_0$ and $\tilde q_0$) are $L^2$
stable with respect to their initial conditions
$Y_0$ and $\tilde Y_0$ in $K$:
\begin{equation}
\label{L2 stability}
||Y(t,\cdot)-\tilde Y(t,\cdot)||_{L^2} 
\le ||Y_{0}-\tilde Y_{0}||_{L^2}+t||q_0-\tilde q_0||_{L^2}.
\end{equation}
for all $t\ge 0$.
This is also true for all $p\ge 1$, when both $Y_0-\tilde Y_0$ and $q_0-\tilde q_0$
belong to $L^p$:
\begin{equation}
\label{Lp stability}
||Y(t,\cdot)-\tilde Y(t,\cdot)||_{L^p} 
\le ||Y_{0}-\tilde Y_{0}||_{L^p}+t||q_0-\tilde q_0||_{L^p}.
\end{equation}
\end{Proposition}

For the sake of completeness, a brief proof of these (standard)
results will be provided at the end of the paper.

\subsection{The main result}

Our main result is
\begin{Theorem}   
\label{main}
Let $Y=Y(t,a,x)$ be a solution to the subdifferential
equation  (\ref{inclusion})
with initial
condition $Y_0\in L^\infty$, with $\partial_a Y_0\ge 0$. Then,
\begin{equation}
\label{solution}  
u(t,y,x)=\int_0^1 H(y-Y(t,a,x))da,
\end{equation}
defines a one parameter family (parameterized by $y\in\rit$)
of Kruzhkov solution to (\ref{scalar}), valued in $[0,1]$. In addition,
all Kruzhkov solutions, with initial values in $L^\infty$, can be
recovered this way (up to a trivial rescaling).
\end{Theorem}

Let us rapidly check the last statement of our main result.
We must show that any Kruzhkov solution $U(t,x)$ with initial
condition $U_0(x)$ valued in $L^\infty$ can
be recovered from a solution to (\ref{inclusion}). 
To do that, according to the first part of the theorem,
it is enough to find an $L^\infty$ 
function $Y_0(a,x)$ such that $\partial_a Y_0\ge 0$ and
$$
U_0(x)=\int_0^1 H(y-Y_0(a,x))da,
$$
for some $y\in\rit$, say $y=1$. This is always possible,
up to rescaling, by assuming:
$$r\le U_0(x)\le 1-r$$
for some constant $r>0$.
Indeed, we set
$$
u_0(y,x)=\max(0,\min(1,y\;U_0(x)))
$$
so that $U_0(x)=u_0(1,x)$ and $\partial_y u_0\ge 0$.
\\
Then, for each fixed $x$,
we solve $u_0(y,x)=a$ by $y=Y_0(a,x)$, setting:
$$
Y_0(a,x)=\frac{a}{U_0(x)},\;\;\;\forall a\in [0,1],\;\;\;\forall x\in\tit^d,
$$
so that 
$$
u_0(y,x)=\int_0^1 H(y-Y_0(x,a))da.
$$
(Notice that $Y_0$ is valued in $[0,r^{-1}]$.)
Finally, according to the first part of the theorem,
we get
$$
U(t,x)=\int_0^1 H(1-Y(t,x,a))da,
$$
where $Y$ is the solution to (\ref{inclusion}) with initial condition $Y_0$.

\subsubsection{Remark}

Notice that, for all $t\ge 0$,
the level sets of $Y$ and $U$ are related by:
$$
\{(a,x),\;\;\;U(t,x)\ge a\}=\{(a,x),\;\;\;Y(t,a,x)\le 1\}.
$$
Thus, the method of construction of $Y_0$ out of $U_0$ and the 
derivation of $U(t,x)$ from $Y(t,a,x)$ can be related
to level-set methods in the spirit of \cite{FSS,Gi1,OF,TGO}. 
This is why we may call
'level-set formulation' of scalar conservation law ({\ref{scalar})
the subdifferential equation given by (\ref{inclusion})

\subsubsection{Remark}

The solutions 
$(t,x)\rightarrow u(t,y,x)$, parameterized by $y\in\rit$,
are automatically ordered in $y$. Indeed, $\partial_y u\ge 0$
immediately follows from representation formula (\ref{solution}).
This is consistent with the order preserving property
of Kruzhkov's theory (as explained in the first section).

\subsection{A second result}

The function $u(t,y,x)$, given by (\ref{solution}),
can also be considered
as a $single$ Kruzhkov solution of a
scalar conservation law in the enlarged
$(1+d)$ dimensional space $\rit\times\tit^d$, namely
\begin{equation}
\label{scalar bis}  
\partial_t u+\partial_y (Q_0(u))+\nabla_x\cdot(Q(u))=0,
\end{equation}
with $(y,x)\in\rit\times\tit^d$,
provided:
\\
\\
1) $Q_0$ is zero,
\\
2) the initial 
condition $u_0(y,x)$ 
is valued in $[0,1]$ and $\partial_y u_0\ge 0$.
\\
Furthermore,
it turns out that, if we add the left-hand side 
$q_0(a)=Q'_0(a)$ to (\ref{inclusion}), so that we get 
(\ref{inclusion bis}):
$$
q_0(a)\;\in\;\;
\partial_t Y+ q(a)\cdot\nabla_x Y+\partial K(Y),
$$
and solve for $Y$, then the corresponding $u$ given by (\ref{solution})
is a Kruzhkov solution to (\ref{scalar bis}).
\\
As a matter of fact, our proof will be done in this larger
framework. We assume that $q_0$, $q$ and $Y_0$ are given in $L^\infty$,
for simplicity. Without loss of generality, up to easy rescalings, we may assume
that both $q_0$ and $Y_0$ are nonnegative, which simplifies some notations.

\begin{Theorem}   
\label{main bis}
Assume that $q_0$ and $q$ are given in $L^\infty$, with $q_0\ge 0$.
Let $Y=Y(t,a,x)$ be a solution to the subdifferential equation
(\ref{inclusion bis}), 
with initial
condition $Y_0\in L^\infty$,
$Y_0\ge 0$ and $\partial_a Y_0\ge 0$. Then,
\begin{equation}
\label{solution bis}  
u(t,y,x)=\int_0^1 H(y-Y(t,a,x))da,
\end{equation}
is the unique Kruzhkov solution to (\ref{scalar bis})
with initial condition:
\begin{equation}
\label{initial bis}  
u_0(y,x)=\int_0^1 H(y-Y_0(a,x))da.
\end{equation}
In addition,
$Y$ is nonnegative and can be recovered from $u$ as:
\begin{equation}
\label{level}
Y(t,a,x)=\int_0^\infty H(u(t,y,x)-a)dy.
\end{equation}
\end{Theorem}

Before proving the theorem, let us observe that the recovery of $Y$ from $u$
through (\ref{level}) is just a consequence of the following elementary lemma
which generalizes (in a standard way) the inversion of a strictly increasing
function of one real variable:
\begin{Lemma}
\label{lemma}
Let: $a\in[0,1]\rightarrow Z(a)\in \rit_+$ with $Z'\ge 0$.
We define the generalized inverse of $Z$:
$$
v(y)=\int_0^1 H(y-Z(a))da,\;\;\;\forall y\in\rit.
$$
Then $v'\ge 0$, $H(y-Z(a))=H(v(y)-a)$ holds true a.e. in $(a,y)\in [0,1]\times\rit$
and:
$$
Z(a)=\int_0^\infty H(a-v(y))dy.
$$
In addtion, for a pair $(Z,v)$, $(\tilde Z,\tilde v)$ of such
functions, we have the co-area formula:
\begin{equation}
\label{coarea}
\int_0^1|Z(a)-\tilde Z(a)|da=
\int_0^1\int_0^\infty|H(y-Z(a))-H(y-\tilde Z(a))|dyda
\end{equation}
$$
=\int_0^1\int_0^\infty|H(v(y)-a)-H(\tilde v(y)-a)|dyda=
\int_0^\infty|v(y)-\tilde v(y)|dy.
$$
\end{Lemma}

To recover (\ref{level}), we notice first that $\partial_a Y\ge 0$ follows from the
very  definition \ref{def} of a solution to (\ref{inclusion bis}). 
Next, $Y\ge 0$ follows from (\ref{maxi})
and the assumptions $q_0\ge 0$, $Y_0\ge 0$. Then, we apply lemma \ref{lemma},
for each fixed $x\in\tit^d$ and $t\ge 0$, by setting $Z(a)=Y(t,a,x)$ and 
$u(t,y,x)=v(y)$.

\subsubsection{Remark}

The function $f(t,a,y,x)=H(y-Y(t,a,x))=H(u(t,y,x)-a)$ valued in $\{0,1\}$
is nothing but the solution
of the Lions-Perthame-Tadmor \cite{LPT}
'kinetic formulation' of (\ref{scalar bis}), which satisfies:
$$
\partial_t f+q_0(a)\partial_y f+q(a)\cdot\nabla_x f=\partial_a \mu,
$$
for some nonnegative measure $\mu(t,a,y,x)$.

\subsubsection{Remark}

As already mentioned, the solutions of (\ref{inclusion bis})
enjoys the $L^p$ stability property with respect to initial
conditions (\ref{Lp stability}), not only for $p=2$ but also 
for all $p\ge 1$.
The case $p=1$ is of particular interest.
Let us consider two solutions $Y$ and $\tilde Y$ of 
(\ref{inclusion bis}) and the corresponding Kruzhkov solutions
$u$ and $\tilde u$ given by Theorem \ref{main bis}.
Using the co-area formula (\ref{coarea}),
we find, for all $t\ge 0$,
$$
\int_{\rit}\int_{\tit^d}|u(t,y,x)-\tilde u(t,y,x)| dxdy=
$$
$$
=\int_0^1\int_{\rit}\int_{\tit^d}|H(u(t,y,x)-a)-H(\tilde u(t,y,x)-a)| dadxdy
$$
$$
=\int_0^1\int_{\rit}\int_{\tit^d}|H(y-Y(t,a,x))-H(y-\tilde Y(t,a,x))| dadxdy
$$
$$
=\int_0^1\int_{\tit^d}|Y(t,a,x)-\tilde Y(t,a,x)| dxda
\le \int_0^1\int_{\tit^d}|Y_0(a,x)-\tilde Y_0(a,x)| dxda
$$
$$
=\int_{\rit}\int_{\tit^d}|u_0(y,x)-\tilde u_0(y,x)| dxdy.
$$

Thus, Kruzhkov's $L^1$ stability property is nothing but a $very$ incomplete
output of the much stronger $L^p$ stability property provided by equation
(\ref{inclusion bis}) for all $p\ge 1$.

\subsubsection{Remark}

As a matter of fact, in Theorem \ref{main bis},
it is possible to translate the $L^p$ stability of
the level set function $Y$ in terms
of the Kruzhkov solution $u$ by using Monge-Kantorovich (MK) distances. 
Let us first recall that
for two probability measures $\mu$ and $\nu$ compactly supported
on $\rit^D$, their $p$ MK distance can be defined (see \cite{Vi}
for instance), for $p\ge 1$, by:
$$
\delta_p^p(\mu,\nu)=\sup\int \phi(x)d\mu(x)+\int \psi(y)d\nu(y),
$$
where the supremum is taken over all pair of continuous functions
$\phi$ and $\psi$ such that:
$$
\phi(x)+\psi(y)\le |x-y|^p,\;\;\;\forall x,y\in\rit^D.
$$
In dimension $D=1$, this definition reduces to:
$$
\delta_p(\mu,\nu)=||Y-Z||_{L^p},
$$
where $Y$ and $Z$ are respectively the generalized inverse (in the sense
of Lemma \ref{lemma}) of 
$u$ and $v$ defined on $\rit$ by:
$$
u(y)=\mu([-\infty,y]),\;\;\;v(y)=\nu([-\infty,y]),\;\;\;\forall y\in\rit.
$$
Next, observe that, for each $x\in\tit^d$, the $y$ derivative of the
Kruzhkov solution $u(t,y,x)$, as described in Theorem \ref{main bis},
can be seen as a probability measure compactly supported on $\rit$.
(Indeed, $\partial_y u\ge 0$, $u=0$ near $y=-\infty$ and $u=1$
near $y=+\infty$.) Then, the $L^p$ stability property simply reads:
$$
\int_{\tit^d} 
\delta_p^p(\partial_y u(t,\cdot,x),\partial_y \tilde u(t,\cdot,x))dx
\le 
\int_{\tit^d} 
\delta_p^p(\partial_y u_0(\cdot,x),\partial_y \tilde u_0(\cdot,x))dx.
$$
We refer to \cite{BBL} and \cite{CFL} for recent occurences
of MK distances in the field of scalar conservation laws.

\section{Proofs}

Let us now prove Theorem \ref{main bis} (which contains the first
part of Theorem \ref{main} as the special case $q_0=0$).
The main idea is to provide, for both formulations (\ref{scalar bis})
and (\ref{inclusion bis}), the same time-discrete approximation scheme, 
namely the 'transport-collapse' method \cite{Br1,Br2,Br3,GM}, 
and get the same limits.

\subsection{A time-discrete approximation}

We fix a time step  $h>0$ and approximate $Y(nh,a,x)$ by $Y_n(a,x)$,
for each positive integer $n$. To get $Y_{n}$ from $Y_{n-1}$, we perform
two steps, making the following induction assumptions:
\begin{equation}
\label{induction}
\partial_a Y_{n-1}\ge 0,\;\;\;0\le Y_{n-1} \le \sup Y_0+(n-1)h\sup q_0,
\end{equation}
which are consistent with our assumptions on $Y_0$.

\subsubsection*{Predictor step}

The first 'predictor' step amounts to solve the linear
equation
\begin{equation}
\label{linear}
\partial_t Y+ q(a)\cdot\nabla_x Y=q_0(a)
\end{equation}
for $nh-h<t<nh$, with $Y_{n-1}$ as initial condition at $t=nh-h$. 
We exactly get at time $t=nh$
the predicted value:
\begin{equation}
\label{predictor}
Y^*_{n}(a,x)=Y_{n-1}(a,x-h\;q(a))+h\;q_0(a).
\end{equation}
Notice that, since $q_0$ is supposed to be nonnegative,
the induction assumption (\ref{induction}) implies:
\begin{equation}
\label{induction bis}
0\le Y^*_{n}\le \sup Y_0+nh\sup q_0.
\end{equation}
However, although $\partial_a Y_{n-1}$ is nonnegative, 
the same may not be true for $\partial_a Y^*_n$. 
This is why, we need a correction step.

\subsubsection*{Rearrangement step}

In the second step, we 'rearrange' $Y^*$ in increasing order with respect to
$a\in [0,1]$, for each fixed $x$, and get the corrected function $Y_{n}$.
Let us recall some elementary facts about rearrangements:

\begin{Lemma}
\label{lemma bis}
Let: $a\in[0,1]\rightarrow X(a)\in \rit_+$ an $L^\infty$
function. Then, there is unique $L^\infty$ function 
$Y:[0,1]\rightarrow \rit_+$, such that $Y'\ge 0$ and:
$$
\int_0^1 H(y-Y(a))da=\int_0^1 H(y-X(a))da,\;\;\;\forall y\in \rit.
$$
We say that $Y$ is the rearrangement of $X$. In addition, for all
$Z\in L^\infty$ such that $Z'\ge 0$, the following
rearrangement inequality:
\begin{equation}
\label{inequality}
\int |Y(a)-Z(a)|^p da\le \int |X(a)-Z(a)|^p da.
\end{equation}
holds true for all $p\ge 1$.
\end{Lemma}

So, we define $Y_n(a,x)$ to be, for each fixed $x$, the rearrangement
of $Y^*_n(a,x)$ in $a\in [0,1]$:
\begin{equation}
\label{corrector}
\partial_a Y_n \ge 0,\;\;\;
\int_0^1 H(y-Y_n(a,x))da=\int_0^1 H(y-Y^*_n(a,x))da,\;\;\;\forall y\in \rit.
\end{equation}
Equivalently, we may define the auxiliary function:
\begin{equation}
\label{u}
u_n(y,x)=\int_0^1 H(y-Y^*_n(a,x))da,\;\;\;\forall y\in \rit,
\end{equation}
i.e.
\begin{equation}
\label{predictor bis}
u_n(y,x)=\int_0^1 H(y-h\;q_0(a)-Y_{n-1}(a,x-h\;q(a)))da,
\end{equation}
and set:
\begin{equation}
\label{corrector bis}
Y_n(a,x)=\int_0^\infty H(a-u_n(y,x))dy.
\end{equation}
At this point, $Y_n$ is entirely determined by $Y_{n-1}$ through formulae
(\ref{predictor}), (\ref{corrector}), or, equivalently, through formulae
(\ref{predictor bis}), (\ref{corrector bis}). 
Notice that, from the very definition 
(\ref{corrector}) of the rearrangement step,
$u_n$, defined by (\ref{u}), can be equivalently written:
\begin{equation}
\label{u bis}
u_n(y,x)=\int_0^1 H(y-Y_n(a,x))da.
\end{equation}
Also notice that, for all function
$Z(a,x)$ such that $\partial_a Z\ge 0$, and all $p\ge 1$:
\begin{equation}
\label{comparison}
\int |Y_n(a,x)-Z(a,x)|^p dadx\le \int |Y^*_n(a,x)-Z(a,x)|^p dadx
\end{equation}
follows from the rearrangement inequality (\ref{inequality}).
Finlly, we see that $\partial_a Y_n \ge 0$ is automatically satisfied
(this was the purpose of the rearrangement step) and
$$
0\le Y_{n}\le \sup Y_0+nh\sup q_0.
$$
follows form (\ref{induction bis}) (since the range of $Y^*_n$ is preserved
by the rearrangement step). So, the induction assumption (\ref{induction})
is enforced at step $n$ and the scheme is well defined.

\subsubsection{Remark}

Observe that, for any fixed $x$, 
$u_n(y,x)$, as a function of $y$, 
is the (generalized) inverse of $Y_n(a,x)$, viewed as a function of $a$,
in the sense of Lemma \ref{lemma}.
Also notice that the level sets $\{(a,y);\;\;y\ge Y_n(a,x)\}$ and 
$\{(a,y);\;\;a\le u_n(y,x)\}$ coincide.

\subsection{The transport-collapse scheme revisited}

The time-discrete scheme can be entirely recast in terms of $u_n$
(defined by (\ref{u bis})). Indeed, introducing
\begin{equation}
\label{tcm1}
ju_n(a,y,x)=H(u_n(y,x)-a),
\end{equation}
we can rewrite (\ref{predictor bis}), (\ref{corrector bis}) in terms
of $u_n$ and $ju_n$ only:
\begin{equation}
\label{tcm2}
u_n(y,x)=\int_0^1 ju_{n-1}(y-h\;q_0(a),x-h\;q(a),a)da.
\end{equation}
We observe that, formulae (\ref{tcm1},\ref{tcm2})
exactly define the 'transport-collapse' (TC) approximation
to (\ref{scalar bis}), or, equivalently, its 'kinetic' approximation,
according to \cite{Br1,Br2,Br3,GM}.

\subsection{Convergence to the Kruzhkov solution}

We are now going to prove that, on one hand, $Y_n(a,x)$ converges to $Y(t,a,x)$
as $nh\rightarrow t$, and, on the other hand, 
$u_n(y,x)$ converges to $u(t,y,x)$,
where $Y$ and $u$ are respectively the unique 
solution to subdifferential equation (\ref{inclusion bis})
with initial condition $Y_0(a,x)$ and the unique Kruzhkov solution to 
(\ref{scalar bis})
with initial condition
\begin{equation}
\label{initial}
u_0(y,x)=\int_0^1 H(y-Y_0(a,x))da.
\end{equation}
\\
From the convergence analysis of the TC method
\cite{Br1,Br2,Br3,GM},
we already know that, as $nh\rightarrow t$, 
$$
\int |u_n(y,x)-u(t,y,x)|dydx\rightarrow 0,
$$
where $u$ is the unique Kruzhkov solution with initial value $u_0$
given by (\ref{initial}). More precisely, if we extend the time
discrete approximations $u_n(y,x)$
to all $t\in [0,T]$ by linear interpolation in time:
\begin{equation}
\label{interpo bis}
u^h(t,y,x)=u_{n+1}(y,x)\frac{t-nh}{h}+u_n(y,x)\frac{nh+h-t}{h},
\end{equation}
then $u^h-u$ converges to $0$ in the space 
$C^0([0,T],L^1(\rit\times\tit^d))$ as $h\rightarrow 0$.
Following (\ref{level}), it is now natural to introduce
the level-set function $Y$ defined by (\ref{level}) from
the Kruzhkov solution:
$$
Y(t,a,x)=\int_0^\infty H(a-u(t,y,x))dy.
$$
(Notice that, at this point, we do not know that $Y$ is a solution to
the subdifferential formulation (\ref{inclusion bis})!)
Let us interpolate the $Y_n$ by
\begin{equation}
\label{interpo}
Y^h(t,a,x)=Y_{n+1}(a,x)\frac{t-nh}{h}+Y_n(a,x)\frac{nh+h-t}{h},
\end{equation}
for all $t\in [nh,nh+h]$ and $n\ge 0$.
By the co-area formula (\ref{coarea}), we have
$$
\int |Y(t,a,x)-Y_n(a,x)|dadx=\int |u(t,y,x)-u_n(y,x)|dydx.
$$
Thus:
$$
\sup_{t\in [0,T]}||Y(t,\cdot)-Y^h(t,\cdot)||_{L^1}
\le \sup_{t\in [0,T]}||u(t,\cdot)-u^h(t,\cdot)||_{L^1}\rightarrow 0,
$$
and we conclude that the approximate solution
$Y^h$ must converge to  $Y$ in $C^0([0,T],L^1([0,1]\times\tit^d))$ 
as $h\rightarrow 0$.
Notice that, since the $Y^h$ are uniformly bounded in $L^\infty$, the convergence
also holds true in $C^0([0,T],L^2([0,1]\times\tit^d))$.

We finally have to prove that $Y$ is the solution to the subdifferential
formulation (\ref{inclusion bis}) with initial condition $Y_0$.

\subsection{Consistency of the transport-collapse scheme}

Let us check that the TC scheme is consistent with
the subdifferential formulation (\ref{inclusion bis}) in its
semi-integral formulation (\ref{semi-integral}).
For each smooth function $Z(t,a,x)$ with $\partial_a Z\ge 0$
and $p\ge 1$, we have
$$
\int |Y_{n+1}(a,x)-Z(nh+h,a,x)|^p dadx 
$$
$$
\le \int |Y_{n+1}^*(a,x)-Z(nh+h,a,x)|^p dadx
$$
(because of property (\ref{comparison})
due to the rearrangement step (\ref{corrector}))
$$
=\int |Y_{n}(a,x-h\;q(a))+h\;q_0(a)-Z(nh+h,a,x)|^p dadx
$$
(by definition of the predictor step (\ref{predictor})
$$
=\int |Y_{n}(a,x)+h\;q_0(a)-Z(nh+h,a,x+h\;q(a))|^p dadx
$$
$$
=\int |Y_{n}-Z(nh,\cdot)|^p dadx+h\;\Gamma+o(h)
$$
where:
$$
\Gamma=p\int (Y_{n}-Z(nh,\cdot))|Y_{n}-Z(nh,\cdot)|^{p-2}
\{q_0-\partial_t Z(nh,\cdot)-q\cdot\nabla_x Z(nh,\cdot)\}dadx
$$
(by Taylor expanding $Z$ about $(nh,a,x)$).
Since the approximate solution provided
by the TC scheme has a unique limit $Y$, as shown in the previous section, 
this limit must satisfy:
$$
\frac{d}{dt} \int |Y-Z|^p dadx
\le
p\;\int (Y-Z)|Y-Z|^{p-2}
(q_0(a)-\partial_t Z-q(a)\cdot\nabla_x Z)dadx,
$$
in the distributional sense in $t$.
In particular, for $p=2$, we exactly recover 
the semi-integral version (\ref{semi-integral}) of (\ref{inclusion bis}).
We conclude that the approximate solutions generated by the TCM
scheme do converge to the solutions of 
(\ref{inclusion bis}) in the sense of Definition \ref{def}, which completes
the proof of Theorem \ref{main bis}.

\section{Viscous approximations}

A natural regularization for subdifferential equation
(\ref{inclusion bis}) amounts to substitute a barrier function for
the convex cone $K$ in $L^2([0,1]\times \tit^d)$ of all functions $Y$
such that $\partial_a Y\ge 0$. Typically, we introduce a convex
function $\phi:\rit\rightarrow ]-\infty,+\infty]$
such that $\phi(\tau)=+\infty$ if $\tau<0$, we define, for all
$Y\in K$,
\begin{equation}
\label{potential}
\Phi(Y)=\int \phi(\partial_a Y)dadx,
\end{equation}
and set $\Phi(Y)=+\infty$ if $Y$ does not belong to $K$.
Typical examples are:
$$\phi(\tau)=-\log(\tau),\;\;\;
\phi(\tau)=\tau\log(\tau),\;\;\;
\phi(\tau)=\frac{1}{\tau},\;\;\;\forall\tau>0.$$
Then, we considered the perturbed subdifferential equation
\begin{equation}
\label{perturbed}
0\;\in\;\;
\partial_t Y+ q(a)\cdot\nabla_x Y-q_0(a)
+\varepsilon\partial \Phi(Y),
\end{equation}
for $\varepsilon>0$.
The general theory of maximal monotone operators guarantees
the convergence of the corresponding solutions to those of
(\ref{inclusion bis}) as $\varepsilon\rightarrow 0$.
It is not difficult (at least formally) to identify
the corresponding perturbation to scalar conservation
(\ref{scalar bis}). Indeed, assuming $\phi(\tau)$ to
be smooth for $\tau>0$, we get, for each smooth function $Y$
such that $\partial_a Y>0$:
$$
\partial\Phi(Y)=-\partial_a(\phi'(\partial_a Y)).
$$
Thus, any smooth solution $Y$ to (\ref{perturbed}),
satisfying $\partial_a Y>0$, solves the
following parabolic equation:
\begin{equation}
\label{viscous}
\partial_t Y+ q(a)\cdot\nabla_x Y-q_0(a)
=\varepsilon\partial_a(\phi'(\partial_a Y)).
\end{equation}
Introducing, the function $u(t,y,x)$ implicitely defined
by 
$$
u(t,Y(t,a,x),x)=a,
$$
we get (by differentiating with respect to $a$, $t$ and $x$):
$$
(\partial_y u)(t,Y(t,a,x),x)\partial_a Y(t,a,x)=1,
$$
$$
(\partial_t u)(t,Y,x)+(\partial_y u)(t,Y,x)\partial_t Y=0,
$$
$$
(\nabla_x u)(t,Y,x)+(\partial_y u)(t,Y,x)\nabla_x Y=0.
$$
Multiplying (\ref{viscous}) by $(\partial_y u)(t,Y(t,a,x),x)$, we get:
\begin{equation}
\label{viscous bis}
-\partial_t u-q(u)\cdot\nabla_x u-q_0(u)\partial_y u
=\varepsilon \partial_y (\phi'(\frac{1}{\partial_y u})).
\end{equation}
In particular, in the case $\phi(\tau)=-\log\tau$,
we recognize a linear viscous approximation to scalar
conservation law (\ref{scalar bis}):
\begin{equation}
\label{viscous ter}
\partial_t u+q(u)\cdot\nabla_x u+q_0(u)\partial_y u
=\varepsilon \partial^2_{yy} u,
\end{equation}
with viscosity only in the $y$ variable. 

\subsubsection{Remark}

Of course, these statements are not rigourous since the parabolic
equations we have considered are degenerate and their solutions
may not be smooth. 

\subsubsection{Remark}

In the case of our main result, Theorem \ref{main}, 
we have $q_0=0$ and the variable $y$ is just a dummy variable 
in (\ref{scalar}).
Thus, the corresponding regularized version
\begin{equation}
\label{viscous quart}
-\partial_t u-q(u)\cdot\nabla_x u
=\varepsilon \partial_y (\phi'(\frac{1}{\partial_y u})).
\end{equation}
includes viscous effects not on the space
variable $x$ but rather on the 'parameter' $y\in\rit$.
This unusual type of regularization has already been used and
analyzed in the level-set framework developped by
Giga for Hamilton-Jacobi equations \cite{Gi2}, and by
Giga, Giga, Osher, Tsai for scalar conservation laws \cite{GG,TGO}.

\section{Related equations}

A similar method can be applied to some special systems of conservation
laws. A typical example (which was crucial for our understanding)
is the 'Born-Infeld-Chaplygin' system considered in \cite{Br4},
and the related concept of 'order-preserving strings'.
This system reads:
\begin{equation}
\label{bi}
\partial_t(hv)+\partial_y(hv^2-hb^2)-\partial_x(hb)=0,
\end{equation}
$$
\partial_t h+\partial_y(hv)=0,\;\;\;
\partial_t (hb)-\partial_x(hv)=0,
$$
where $h,b,v$ are real valued functions of time $t$ and two space
variables $x,y$.
In \cite{Br4}, this system is related to the following
subdifferential system:
\begin{equation}
\label{bi-subdif}
0\in \partial_t Y-\partial_x W+\partial K(Y),
\;\;\;\partial_t W=\partial_x Y,
\end{equation}
where $(Y,W)$ are real valued functions of $(t,a,x)$ 
and $K$ is the convex
cone of all $Y$ such that $\partial_a Y\ge 0$.
The (formal) correspondence between (\ref{bi}) and (\ref{bi-subdif})
is obtained by setting:
$$
h(t,x,Y(t,x,a))\partial_a Y(t,x,a)=1,
$$
$$
v(t,x,Y(t,x,a))=\partial_t Y(t,x,a),\;\;\;
b(t,x,Y(t,x,a))=\partial_x Y(t,x,a).
$$
Unfortunately, this system is very special (its smooth solutions
are easily integrable). In our opinion, it is very unlikely that 
$L^2$ formulations can be found for general
hyperbolic conservation laws as easily as in the multidimensional
scalar case.

\section{Appendix: proof of Proposition \ref{Proposition}}

In the case when $q_0$ and $Y_0$ belong
to $L^\infty$ and are nonnegative,
we already know, from the convergence of
the TC scheme, that there is a solution $Y$ to (\ref{inclusion bis}),
with initial value $Y_0$, in the sense of definition \ref{def}.
From (\ref{induction}), we also get for such solutions, 
when $q_0\ge 0$ and  $Y_0\ge 0$,
$$
0\le Y(t,\cdot)\le \sup Y_0+t\sup q_0,\;\;\;\forall t\ge 0.
$$
By elementary rescalings, we can remove the assumptions that 
both $Y_0$ and $q_0$ are nonnegative and get estimate (\ref{maxi}).
\\
Let us now examine some additional properties of the solutions to
(\ref{inclusion bis}) obtained from the TC approximations.
First, we observe that, in the TC scheme,
\\
1) the predictor step (a translation in the $x$ variable by
$h\;q(a)$ plus an addition of $h\;q_0(a)$) is isometric in all $L^p$ spaces,
\\
2) the corrector step (an increasing rearrangement in the $a$ variable)
is non-expansive in all $L^p$.
\\
Thus the scheme is non-expansive in all $L^p([0,1]\times \tit^d)$.
More precisely, for two different initial conditions $Y_0$ and $\tilde Y_0$,
and two different data $q_0$ and $\tilde q_0$, all in $L^\infty$,
we get for the corresponding approximate solutions $Y_n$ and $\tilde Y_n$:
\begin{equation}
\label{non expansive}
||Y_{n}-\tilde Y_{n}||_{L^p} 
\le ||Y_{n-1}-\tilde Y_{n-1}||_{L^p}+h||q_0-\tilde q_0||_{L^p}\;.
\end{equation}
This shows that (\ref{Lp stability}) holds true for all solutions
of (\ref{inclusion bis}) generated by the TC scheme.
\\
Since the scheme is also invariant under translations in the $x$ variable,
we get the following a priori estimate:
\begin{equation}
\label{esti2}
||\nabla_x Y_n||_{L^p}\le ||\nabla_x Y_0||_{L^p}.
\end{equation}
Finally, let us compare two solutions of the scheme $Y_n$ and
$\tilde Y_n=Y_{n+1}$ obtained with initial condition $\tilde Y_0=Y_1$.
Using (\ref{non expansive}), we deduce:
$$
\int |Y_{n+1}(a,x)- Y_{n}(a,x)|^p dadx
\le \int |Y_{1}(a,x)- Y_{0}(a,x)|^p dadx
$$
$$
\le \int |Y^*_{1}(a,x)- Y_{0}(a,x)|^p dadx
=\int |Y_{0}(a,x-h\;q(a))+h\;q_0(a)-Y_0(a,x)|^p dadx.
$$
So we get a second a priori estimate:
\begin{equation}
\label{esti3}
||Y_{n+1}-Y_n||_{L^p}\le 
(||q_0||_{L^p}+||q||_{L^\infty}||\nabla_x Y_0||_{L^p})h.
\end{equation}
Thus the solutions $Y$ to (\ref{inclusion bis}) obtained from the TC
scheme satisfy the a priori bounds:
\begin{equation}
\label{esti2 bis}
||\nabla_x Y(t,\cdot)||_{L^p}\le ||\nabla_x Y_0||_{L^p},
\end{equation}
\begin{equation}
\label{esti3 bis}
||\partial_t Y(t,\cdot)||_{L^p}\le 
||q_0||_{L^p}+||q||_{L^\infty}||\nabla_x Y_0||_{L^p}.
\end{equation}
Notice that, at this level, we still do not know if solutions,
in the sense of Definition \ref{def} exist when
$Y_0\in K$ and $q_0\in L^2([0,1])$ are not in $L^\infty$ and we know nothing about their
uniqueness. This can be easily addressed by standard functional analysis
arguments. 

\subsubsection*{Existence for general data}

Let $Y_0\in K$ and $q_0\in L^2([0,1])$. We can find two Cauchy sequences in $L^2$, labelled
by $k\in\nit$, namely
$Y_0^k\in K$ and $q_0^k\in L^2([0,1])$, made of smooth functions, with limits $Y_0$ and $q_0$
respectively. Let us denote by $Y^k$ the corresponding solutions, generated by
the TC scheme. Because of their $L^2$ stability, they satisfy:
$$
\sup_{t\in [0,T]}||Y^k(t,\cdot)-Y^{k'}(t,\cdot)||_{L^2}
\le ||Y_0^k-Y_0^{k'}||_{L^p}+T||q_0^k-q_0^{k'}||_{L^2}.
$$
So, $Y^k$ is a Cauchy sequence in $C^0([0,T],L^2)$ of solutions of
(\ref{inclusion bis}) in the sense of Definition \ref{def}, with a definite limit $Y$.
Definition \ref{def} is clearly stable under this convergence process.
So, we conclude that $Y$ satisfies the requirements of Definition \ref{def}
and is a solution with initial condition $Y_0$ and left-hand side $q_0$.
Notice that, through our approximation process, we keep the a priori estimates
(\ref{esti2 bis}),(\ref{esti3 bis}), for general data $q_0\in L^2([0,1])$.

\subsubsection*{Uniqueness}

Let us consider a solution $Y$ to (\ref{inclusion bis}),
with initial condition $Y_0\in K$ and
left-hand side $q_0\in L^2([0,1])$, in the sense of Definition \ref{def}.
By definition $Y(t,\cdot)\in K$ depends continuously of $t\in [0,T]$ in $L^2$.
From definition (\ref{semi-integral}), using $Z=0$ as a test function,
we see that:
$$
\frac{d}{dt}||Y(t,\cdot)||^2_{L^2} \le 2\int Y(t,a,x)q_0(a)\; dadx \le 
||Y(t,\cdot)||^2_{L^2}+||q||^2_{L^2},
$$
which implies that the $L^2$ norm $Y(t,\cdot)$ stays uniformly bounded on any finite
interval $[0,T]$. Thus, $T>0$ being fixed,
we can mollify $Y$ and get, for each $\epsilon\in ]0,1]$
a smooth function $Y_\epsilon$, valued in $K$,
so that:
\begin{equation}
\label{error}
\sup_{t\in [0,T]}||Y(t,\cdot)-Y_\epsilon(t,\cdot)||_{L^2}\le \epsilon.
\end{equation}
Let us now consider an initial condition $Z_0$ such that $\nabla_x Z_0$
belongs to $L^2$. We know that there exist a solution $Z$ to 
(\ref{inclusion bis}), still in the sense of Definition \ref{def},
obtained by TC approximation,
for which both $\partial_t Z(t,\cdot)$ and $\nabla_x Z(t,\cdot)$ 
stay uniformly bounded in $L^2$ for all $t\in [0,T]$. This function
$Z$ has enough regularity to be used as 
a test function in (\ref{semi-integral}) when expressing that $Y$ is a solution
in the sense of Definition \ref{def}.
So, for 
each smooth nonnegative function $\theta(t)$, compactly supported in $]0,T[$,
we get from (\ref{semi-integral}):
$$
\int \{\theta'(t)|Y-Z|^2
+2\theta(t)(Y-Z)(q_0(a)-\partial_t Z- q(a)\cdot\nabla_x Z)\}dadxdt\ge 0.
$$
Substituting $Y_\epsilon$ for $Y$, we have, thanks to estimate (\ref{error}),
$$
\int \{\theta'(t)|Y_\epsilon-Z|^2
+2\theta(t)(Y_\epsilon-Z)(q_0(a)-\partial_t Z- q(a)\cdot\nabla_x Z)\}dadxdt
\ge -C\epsilon,
$$
where $C$ is a constant depending on $\theta$, $Z$, $q_0$ and $q$ only.
Since $Z$ is also a solution, using $Y_\epsilon$ as a test function, we get
from formulation (\ref{semi-integral}):
$$
\int \{\theta'(t)|Z-Y_\epsilon|^2
+2\theta(t)(Z-Y_\epsilon)
(q_0(a)-\partial_t Y_\epsilon- q(a)\cdot\nabla_x Y_\epsilon)\}dadxdt
\ge 0.
$$
Adding up these two inequalities, we deduce:
$$
\int\{2\theta'(t)|Y_\epsilon-Z|^2
+2\theta(t)(Y_\epsilon-Z)(\partial_t (Y_\epsilon-Z)+ q(a)\cdot\nabla_x(Y_\epsilon- Z))\}dadxdt
\ge -C\epsilon.
$$
Integrating by part in $t\in [0,T]$ and $x\in\tit^d$, we simply get:
$$
\int \theta'(t)|Y_\epsilon-Z|^2 dadxdt
\ge -C\epsilon.
$$
Letting $\epsilon\rightarrow 0$, we deduce:
$$
\frac{d}{dt}\int |Y-Z|^2 dadx\le 0.
$$
We conclude, at this point, that:
$$
||Y(t,\cdot)-Z(t,\cdot)||_{L^2}\le ||Y_0-Z_0||_{L^2},\;\;\;\forall t\in [0,T]
$$
This immediately implies the uniqueness of $Y$. Indeed, any other solution $\tilde Y$
with initial condition $Y_0$ must also satisfy:
$$
||\tilde Y(t,\cdot)-Z(t,\cdot)||_{L^2}\le ||Y_0-Z_0||_{L^2}.
$$
Thus, by the triangle inequality:
$$
||\tilde Y(t,\cdot)-Y(t,\cdot)||_{L^2}\le 2||Y_0-Z_0||_{L^2}.
$$
Since $Z_0\in K$ is any function such that $\nabla_x Z_0$ belongs to $L^2$,
we can make $||Y_0-Z_0||_{L^2}$ arbitrarily small and conclude that $\tilde Y=Y$,
which completes the proof of uniqueness.

\subsection*{Acknowledgments}
This article was written at the Bernoulli Centre, 
EPFL, Lausanne, in September 2006,
during the program ``Asymptotic Behaviour in Fluid Mechanics''.
The author is grateful to the organizers, Dragos Iftime, 
Genevi\`eve Raugel and Tudor Ratiu for their kind invitation.


\begin{thebibliography}{B}
\bibitem[BB]{BB} S. Bianchini, A. Bressan, 
{\sl Vanishing viscosity solutions of nonlinear hyperbolic systems,}
{\it  Ann. of Math. (2)  161  (2005) 223-342.}


\bibitem[BDLL]{BDLL} G. Boillat. C. Dafermos, P. Lax, T.P. Liu,
{\sl Recent mathematical methods in nonlinear wave propagation,}
{\it Lecture Notes in Math., 1640, Springer, Berlin, 1996}


\bibitem[BBL]{BBL} F. Bolley, Y. Brenier, G. Loeper,
{\sl Contractive metrics for scalar conservation laws,}
{\it J. Hyperbolic Differ. Equ. 2 (2005) 91-107.}

\bibitem[Br1]{Br1} Y. Brenier,
{\sl Une application de la sym\'etrisation de Steiner aux \'equations 
hyperboliques: la m\'ethode de transport et \'ecroulement,}
{\it  C. R. Acad. Sci. Paris Ser. I Math.  292  (1981) 563-566.}

\bibitem[Br2]{Br2} Y. Brenier,
{\sl R\'esolution d'\'equations d'\'evolution quasilin\'eaires en dimension $N$
 d'espace \`a l'aide d'\'equations lin\'eaires en dimension $N+1$,}
 {\it J. Differential Equations  50  (1983) 375-390.}

 \bibitem[Br3]{Br3} Y. Brenier,
 {\sl Averaged multivalued solutions for scalar conservation laws,}
 {\it  SIAM J. Numer. Anal.  21  (1984) 1013-1037.}


 \bibitem[Br4]{Br4} Y. Brenier,
 {\sl Order preserving vibrating strings and applications to electrodynamics 
 and magnetohydrodynamics,}
 {\it  Methods Appl. Anal.  11  (2004) 515-532.}

 \bibitem[Brn]{Brn} P. Brenner,
 {\sl The Cauchy problem for symmetric hyperbolic systems in $L\sb{p}$,}
 {\it Math. Scand. 19 (1966) 27-37.}

 \bibitem[Brz]{Brz} H. Brezis,
 {\it Op\'erateurs maximaux monotones et semi-groupes de contractions 
 dans les espaces de Hilbert,}
 {\sl North-Holland Mathematics Studies, No. 5. 
 1973.}

\bibitem[CFL]{CFL} 
J. A. Carrillo, M. Di Francesco, C. Lattanzio,
{\sl Contractivity of Wasserstein Metrics and Asymptotic Profiles for Scalar 
Conservation Laws,}
{\it Preprints on Conservation Laws, 2006,
http://www.math.ntnu.no/conservation/2006/003.html.}

\bibitem[Da]{Da} C. Dafermos,
{\sl Hyperbolic conservation laws in continuum physics,}
{\it Springer-Verlag, Berlin, 2000.}

\bibitem[Di]{Di} R. DiPerna,
{\sl Convergence of approximate solutions to conservation laws,}
{\it Arch. Rational Mech. Anal. 82 (1983) 27-70.}




\bibitem[FSS]{FSS} R. Fedkiw, G. Sapiro, C-W Shu,
{\it shock capturing, level sets, and PDE based methods in computer vision 
and image processing: a review of Osher's contributions,}
{\sl  J. Comput. Phys.  185  (2003) 309-341.}

\bibitem[GG]{GG} M-H. Giga, Y. Giga,
{\it Minimal vertical singular diffusion preventing 
overturning for the Burgers equation,}
{\sl  Recent advances in scientific computing and PDEs,
Contemp. Math., 330, Amer. Math. Soc., 2003.}

\bibitem[Gi1]{Gi1} Y. Giga,
{\sl Surface evolution equations. A level set approach,}
{\it Monographs in Mathematics, 99. Birkhäuser Verlag, Basel, 2006.}

\bibitem[Gi2]{Gi2} Y. Giga,
{\it Viscosity solutions with shocks,}
{\sl Comm. Pure Appl. Math.  55  (2002) 431-480.}

\bibitem[GM]{GM} Y. Giga, T. Miyakawa, 
 {\sl A kinetic construction of global solutions of first order quasilinear 
 equations,}
 {\it Duke Math. J.  50  (1983) 505-515.}


\bibitem[Kl]{Kl} S. Klainerman, 
{\it The null condition and global existence to nonlinear wave equations,}
{\sl Nonlinear systems of partial differential equations in applied 
mathematics, 293-326, Lectures in Appl. Math., 23, Amer. Math. Soc., 1986.}


\bibitem[Kr]{Kr} S. N. Kruzhkov,
{\it First order quasilinear equations with several independent variables,}
{\sl Mat. Sb. (N.S.) 81 (123) (1970) 228-255.}

\bibitem[LPS]{LPS} P.-L. Lions, B. Perthame, T. Souganidis, 
{\sl Existence and stability of entropy solutions for the hyperbolic systems 
of isentropic gas dynamics in Eulerian and Lagrangian coordinates,}
{\it  Comm. Pure Appl. Math.  49  (1996) 599-638.}


\bibitem[LPT]{LPT} P.-L. Lions, B. Perthame, E. Tadmor,
{\sl A kinetic formulation of multidimensional scalar conservation laws and 
related equations,}
{\it  J. Amer. Math. Soc.  7  (1994) 169-191.}

\bibitem[Ma]{Ma} A. Majda, 
{\sl Compressible fluid flow and systems of conservation laws in several 
space variables,}
{\it Applied Mathematical Sciences, 53. Springer-Verlag, 1984.}

\bibitem[OF]{OF} S. Osher, R. Fedkiw,
{\it Level set methods.  Geometric level set methods in imaging, 
vision, and graphics,}
{\it Springer, New York, 2003.}




 \bibitem[Se]{Se} D. Serre,
 {\it Systems of conservation laws,}
 {\sl Cambridge University Press, Cambridge, 2000.}


\bibitem[Ta]{Ta} L. Tartar, 
{\sl Compacit\'e par compensation: r\'esultats et perspectives,}
{\sl Nonlinear partial differential equations and their applications,
Res. Notes in Math., 84, Pitman, Boston 1983.}


 \bibitem[TGO]{TGO} Y-H. R. Tsai, Y. Giga, S. Osher, 
 {\sl A level set approach for computing discontinuous solutions of 
 Hamilton-Jacobi equations,}
 {\it  Math. Comp.  72  (2003) 159-181.}

\bibitem[Vi]{Vi} C. Villani, 
{\sl Topics in optimal transportation,}
{\it American Mathematical Society, Providence, 2003.}

 \end{thebibliography}
 \end{document}